\documentclass[11pt]{amsart}

\usepackage{amsmath,amsmath,amssymb,latexsym}

\newtheorem{theorem}{Theorem}
\newtheorem{lemma}{Lemma}
\newtheorem{corollary}{Corollary}

\newtheorem{proposition}{Proposition}
\newtheorem{result}[proposition]{Result}

	\usepackage{setspace}
	\makeatletter
	\DeclareRobustCommand*\cal{\@fontswitch\relax\mathcal}
	\makeatother

	\usepackage{times, amsmath}
	\usepackage{amsfonts}
	\usepackage{amssymb}
	\usepackage{amsthm}
        \usepackage{wrapfig}
	\usepackage[usenames]{color}
	\usepackage[table,xcdraw]{xcolor}
        \usepackage{graphicx}
        \usepackage{multicol}
        \usepackage{multirow}
        \usepackage{natbib}
        \usepackage{url}
     \usepackage{booktabs}
\usepackage[T1]{fontenc}
\usepackage{times}

\usepackage{endnotes}
\usepackage{enumerate}
\usepackage{subfigure}
\usepackage{wrapfig}

\newcommand{\ba}{\begin{array}}
\newcommand{\ea}{\end{array}}
\newcommand{\be}{\begin{equation}}
\newcommand{\ee}{\end{equation}}
\newcommand{\bea}{\begin{eqnarray}}
\newcommand{\eea}{\end{eqnarray}}

\newcommand{\chg}{\color{black}}

\DeclareMathOperator*{\sgn}{sgn}
\DeclareMathOperator*{\Epi}{Epi}

\begin{document}

\title{Penalized Euclidean Distance Regression}

\author{D. VASILIU}
\address{Department of Mathematics and Applied Mathematics, Virginia Commonwealth University, Richmond, Virginia, 23284, USA} 
\email{dvasiliu@vcu.edu}

\author{T. DEY}
\address{Department of Quantitative Health Sciences, Lerner Research Institute, Cleveland Clinic, Cleveland, Ohio, 44195, USA} 
\email{deyt@ccf.org}

\author{I.L. DRYDEN}
\address{School of Mathematical Sciences, The University of Nottingham, Nottingham, NG7 2RD, UK}\email{ian.dryden@nottingham.ac.uk}

\date{}

\maketitle


\begin{abstract}
A new method is proposed for variable selection and prediction in 
linear regression problems where the 
number of predictors can be much larger than the number of observations. The method involves minimizing 
a penalized Euclidean distance (equivalent in many properties to the empirical norm), where the penalty is the geometric mean of the $\ell_1$ and $\ell_2$ 
norms of the regression coefficients. This particular formulation exhibits a  grouping effect, 
which is useful for model selection in high dimensional problems.  
Also, an important result is a signal recovery theorem, which does not require an estimate of the 
noise standard deviation.  
Practical performances of variable selection and prediction are evaluated through simulation studies and 
the analysis of a couple of real datasets.  
\end{abstract}

\keywords{
Euclidean distance; Grouping; Penalization; Prediction; Regularization; Sparsity; Variable screening.
}

\section{Introduction}
High dimensional regression problems are of great interest in a wide range of applications, 
for example in analysing microarrays \citep{Hastieetal:08,Fanetal:09}, functional magnetic resonance images 
\citep{Gaudesetal:13} and mass spectrometry data \citep{Tibshiranietal:05}. 
We consider the problem of predicting a single response $Y$ from a set of $p$ predictors $X_1,\ldots,X_p$, 
where $p$ can be much larger than the number of observations $n$ of each variable. If $p > n$, commonly used methods include regularization by adding a penalty to the least squares objective function or variable
selection of the most important predictors. 

A wide range of methods is available 
for achieving one or both of the essential goals in linear regression: accomplishing predictive accuracy and identifying
pertinent predictive variables. 
There is a very large literature on high-dimensional regression methods, for example introductions to the area are 
given by \cite[Section 3.4]{Hastieetal:08} and \cite[Chapter 3]{Jamesetal:13}. 
Earlier methods for high-dimensional regression include procedures which minimize a least squares objective function plus a penalty on the 
regression parameters. The methods include 
ridge regression \citep{Hoerl:70a,Hoerl:70b} with a squared $\ell_2$ penalty; Lasso \citep{Tibshirani:96} with an $\ell_1$ penalty; 
and the Elastic Net \citep{Zou:05} with a linear combination of 
$\ell_1$ and squared $\ell_2$ penalties. Alternative methods include the Dantzig selector \citep{Candes:07}, where the correlation between the residuals 
and predictors is bounded, Sure Independence Screening \citep{Fan:08} where predictors are initially screened 
using componentwise regression, and Square Root Lasso \citep{Belloni}, which involves minimizing an empirical norm of the residuals 
 with an $\ell_1$ penalty. 

In our method we use the Euclidean distance objective function 
plus a new norm based on the geometric mean of the $\ell_1$  and $\ell_2$ norms of the regression parameters. The advantage of our approach is that we are also able to provide  
the pivotal recovery property, but in addition gain the grouping property of the Elastic Net (where regression
coefficients of a group of highly correlated variables are very similar). The resulting penalized Euclidean distance method is shown to work well in a variety of settings. A particularly strong feature 
is that it works well when there are correlated designs with 
weak signal and strong noise. 

\section{Penalized Euclidean distance} \label{PED}
\subsection{Notation}
We assume that the data are organized as an $n\times p$ design matrix $X$, and a $n$ dimensional response vector $Y$, where $n$ is the number of observations and $p$ is the number of variables. The columns of the matrix $X$ are denoted by $X_j$, i.e. $X_j=(x_{1,j},x_{2,j}...,x_{n,j})^T$, $j=1,...,p$ and the regression parameters are $\beta = (\beta_1,\ldots,\beta_p)^T$. 
We assume that a vector of outcomes $Y$ is modelled linearly as 
$
Y = X \beta^* + \sigma \epsilon
$
where $\beta^*$ is the true parameter vector of dimension $p$, the expectation of $\epsilon$ is zero and its variance is the identity matrix.
Thus we assume that the expectation of the response  $Y=(y_1,\ldots,y_n)^T$ depends only on a few variables, and so
\begin{equation}
X \beta^* =  X^* \tilde\beta^* ,  \label{drivers}
\end{equation} 
where the columns of the matrix $X^*$ are a subset of the set of columns of the entire design matrix $X$, so $X^*$ is associated with a subset of indices ${\mathfrak{J}^*}\subset\{1,2,\ldots,p\}$ and $\tilde\beta^*$ is a subvector of $\beta^*$ with the zero elements removed whose dimension is equal to the cardinality of $\mathfrak{J}^*.$  
In general, if we try to minimize  
$\| Y - X^* \tilde\beta^* \|$ over choices of $\mathfrak{J}^*$ and vectors
$\tilde\beta^*$, the optimal choice of  
$\mathfrak{J}^*$ may not be unique since an under-determined system could have solutions with different sparsity patterns, even if the degree of the optimal sparsity (model size) is the same. However, in the signal reconstruction problem that we consider where a penalty on 
the parameters is introduced, 
we will show that under some assumptions we can approximate $\beta^*$ in probability. The cardinality of $\mathfrak{J}^*$ (denoted by $|\mathfrak{J}^*|$) is assumed to be less than the number of observations and when $p$ is much greater than $|\mathfrak{J}^*|$ a huge challenge is to detect the set of irrelevant columns, namely the variables that correspond to the position of the null components of $\beta^*$ and thus, not needed for efficiently controlling the outcomes $Y$. 

\subsection{The Penalised Euclidean Distance objective function}
Our method involves minimizing the Euclidean distance (as a loss function, essentially equivalent to the empirical norm used in {\em Square Root Lasso}) between 
$Y$ and $X \beta$, with a penalty based on the geometric 
mean of the $\ell_1$ and $\ell_2$ norms. In particular, 
we minimize
\begin{equation}
L_{PED}(\lambda,\beta)=\|Y-X\beta\|+\lambda\sqrt{\|\beta\|  \|\beta\|_1}   \label{obj}
\end{equation}
where $\lambda$ is scalar regularization parameter, $\beta=(\beta_1,\beta_2...,\beta_p)$ is a vector in $\mathbb{R}^p$ 
(to be optimized over), $\|\beta\|^2=\sum\limits_{i=1}^{p}\beta^2_j$ is the squared $\ell_2$ norm and $\|\beta\|_1=\sum\limits_{i=1}^{p}|\beta_j|$, the $\ell_1$ norm. The Penalized Euclidean Distance estimator $\hat\beta$ is defined as the minimizer of the objective function (\ref{obj}), i.e. $\hat\beta=(\hat\beta_1,\hat\beta_2,\ldots,\hat\beta_p)$ and
\begin{equation}\label{optim}
\hat\beta(\lambda)=\arg\min\limits_{\beta \in \mathbb{R}^p}\{L_{PED}(\lambda,\beta)\}.
\end{equation}
The penalty is proportional to the geometric mean of the $\ell_1$ and $\ell_2$ norms and has only one control parameter, $\lambda$.

An alternative, well-established method that involves a convex combination $\ell_1$ and $\ell_2^2$ penalties is the {elastic net} \citep{Zou:05}. This is based on the {na\"{i}ve elastic net} criterion, whose objective function is defined as
\begin{equation*}
L_{nen}(\lambda_1,\lambda_2,\beta)=\|Y-X\beta\|^2+\lambda_2\|\beta\|^2+\lambda_1\|\beta\|_1
\end{equation*}
where $\hat\beta_{en}=\sqrt{1+\lambda_2}\hat\beta_{nen}$ and $\hat\beta_{nen}=\arg\min\limits_{\beta}\{L_{nen}(\lambda_1,\lambda_2,\beta)\}$. The Lasso \citep{Tibshirani:96} is
a special case with $\lambda_1 > 0, \lambda_2 = 0$ and ridge regression has $\lambda_1 = 0, \lambda_2 > 0$, and so the {\em elastic net} combines the two methods. 
Our method also combines features of Lasso and Ridge Regression but in a radically different way.  
The Square Root Lasso \citep{Belloni} involves minimizing 
\begin{equation}
L_{SQL}(\lambda,\beta)=\frac{1}{n} \|Y-X\beta\|+\frac{\lambda}{n} \|\beta\|_1   \label{SQLobj}
\end{equation} 
and so the first term for the Penalized Euclidean Distance estimator is the same as that of Square Root Lasso multiplied by $n$, and it is just the penalty that differs. \citet{Belloni} 
have given a rationale for choosing the regularization parameter using a property called pivotal recovery, 
without requiring an estimate of the 
noise standard deviation.  

The Penalized Euclidean Distance penalty is identical to the Lasso penalty for a single non-zero $\beta_i$, and so for very sparse 
models behaviour like the Square Root Lasso is envisaged. 

We shall show that for the Penalized Euclidean Distance estimator there is a grouping effect for correlated variables, which is a property shared by the {\em elastic net}. A grouping effect occurs where highly correlated predictors $X_j, X_k$ will give rise to very similar regression 
parameter estimates, i.e. $\hat\beta_j \approx \hat\beta_k$.

\subsection{Standardizing to the unit hyper-sphere}
By applying a location transformation, both the design matrix $X$ and the response vector $Y$ can be centred, and we also scale the predictors so that
\begin{equation}
\sum\limits_{i=1}^{n}y_i=0,\,\,\,\sum\limits_{i=1}^{n}x_{i,j}=0,\,\,\,\sum\limits_{i=1}^{n}x^2_{i,j}=1,\,\,\,j=1,...,p. \label{stand}
\end{equation}
Each covariate $X_j$ can be regarded as a point on the unit hypersphere $S^{n-1}$ with a centring constraint. 
We assume that the global minimum of $\|Y-X^*\beta\|$ is very small, but nonetheless positive 
\begin{equation*}
\min\limits_{\beta\in\mathbb{R}^{|\mathfrak{J}^*|}} \|Y - X^*\beta\|\geq c' >0, 
\end{equation*}
as obtained in the presence of noise and when the number of observations is larger than $|\mathfrak{J}^*|$, and $ X^*, |\mathfrak{J}^*|$ were defined after (\ref{drivers}).

For any vector $\beta\in\mathbb{R}^p$, we denote by $\theta_j$ the angle between $X_j$ and $Y-X\beta$. Thus for a vector $\beta$ that achieves the global minimum of $\|Y - X^*\beta\|$, 
we must have $\theta_j=\frac{\pi}{2}$ for $j \in \mathfrak{J}^*$ and we can define the set of solutions as
\begin{equation*} 
\mathcal{S}=\{\beta\in\mathbb{R}^p\,:\, \theta_j=\frac{\pi}{2}\,\, \text{for all}\, j\in\mathfrak{J}^*\,\,\text{and}\,\beta_j=0\,\, \text{for all}\, j\in{\mathfrak{J}^*}^c\}
\end{equation*}
where ${\mathfrak{J}^*}^c$ denotes the complement of $\mathfrak{J}^*$ in the set of all indexes $\{1,2,\ldots,p\}.$
\noindent We also assume that $\mathcal{S}$ is bounded away from $0_{\mathbb{R}^p}$ (i.e. $\beta=0_{\mathbb{R}^p}$ does not minimize $\|Y-X\beta\|$).  Our goal is to build an estimator of the index set $\mathcal{S}.$ Also, thinking of covariates as vectors on the unit hypersphere $S^{n-1}$, we would like to build an objective function that facilitates automatic detection of overcrowding, which is a situation where there is a group of very close covariates on the hypersphere
(where the great circle distances are small within the group) and these correspond to highly correlated predictors. 
If an estimation method exhibits a grouping effect then in the 
predictors in this group  will have similar estimated regression parameters.  
Throughout the paper we assume that $Y$ have been centred and the columns of $X$ have been standardised as described above.

\section{Theoretical Results}\label{Theory}
\subsection{Geometric mean norm and grouping}
The Penalized Euclidean Distance objective function enables variable selection under some mild compatibility conditions. 
The concept is based on the simple fact that the sum of the squares of the relative sizes of vector components (as defined by $\beta_j/\|\beta\|_2$)  is always equal to 1. For any vector in $\mathbb{R}^p$, if there are
components that have relative size larger than $\frac{1}{\sqrt{p}}$ then the other components must have relative size falling under this value. In addition if many components 
have similar relative size due to a grouping effect, then the relative size of those components must be small. 
The new penalty function that we consider is actually a norm. 

\begin{lemma}\label{convlema}
Given any two $p$-norms $f_{p_1},f_{p_2}:\mathbb{R}^n\rightarrow [0,+\infty)$, i.e., for some $p_1,p_2\geq 1$, $f_{p_1}(\beta)=\left( \displaystyle\sum\limits_{i=1}^{n}|\beta_i|^{p_1}\right)^{1/p_1}$, $f_{p_2}(\beta)=\left( \displaystyle\sum\limits_{i=1}^{n}|\beta_i|^{p_2}\right)^{1/p_2}$, we have that $\sqrt{f_{p_1}  f_{p_2}}$ is also a norm on $\mathbb{R}^n$.
\end{lemma}

\noindent The following theorem demonstrates the grouping effect achieved by a minimizer of the penalized Euclidean distance objective function. The idea of grouping effect was first introduced by \citet{Zou:05}. Our version of the grouping effect involves the relative contributions of the components of the minimizer of the Penalized Euclidean Distance objective function. This property enables the process of eliminating irrelevant variables from the model.

Considering the situation of very large $p$ compared to $n$, selecting and grouping variables is an important priority. Theorem \ref{grouping} below supports the idea of obtaining groups of highly correlated variables, based on the relative size of the corresponding component minimizers of the penalized Euclidean distance objective function. 

\begin{theorem}\label{grouping}
Assume we have a standardized data matrix $X$, and $Y$ is a centred response vector, 
as in (\ref{stand}). Let $\hat\beta$ be the Penalized Euclidean Distance estimate given by 
$$\hat\beta(\lambda)=\arg\min\limits_{\beta}\{L_{PED}(\lambda,\beta)\}$$ 
for some $\lambda>0$. Define
\[
D_{\lambda}(i,j)=\frac{1}{\|\hat\beta(\lambda)\|}|\hat\beta_i(\lambda)-\hat\beta_j(\lambda)|
\]
then
\[
D_{\lambda}(i,j)\leq \frac{{\chg 2}(1-\rho_{ij})^{1/2}}{\lambda} \leq\frac{{\chg 2}\theta_{ij}}{\lambda}
\]
where $\rho_{ij} = (X_i)^T(X_j)$,  is the sample correlation, $\theta_{ij}$ is the angle between $X_i$ and $X_j$, {\chg $0 \le \theta_{ij} \le \pi/2$}.
\end{theorem}

\noindent Note that this result is analogous to Theorem 1 of \citet{Zou:05} for the Elastic Net, and the same method of
proof is used in the Appendix. 

From Theorem \ref{grouping} if $\theta_{ij}$ is small then the corresponding parameters estimated from Penalized Euclidean Distance regression will be similar, which is the grouping effect.  
When $\theta_{ij}=0$ we have the Corollary:

\begin{corollary}\label{group}
Let $\hat\beta(\lambda)=\arg\min\limits_{\beta}\{L_{PED}(\lambda,\beta)\}.$  If $X_i = X_j$ then $\hat\beta_i(\lambda)=\hat\beta_j(\lambda).$
\end{corollary}

The grouping effect means that if we have 
strong overcrowding on the unit hypersphere around an irrelevant column then this would be detected by a large drop in 
the relative size of the corresponding components of the solution to our objective function.

\subsection{Sparsity Properties}
We consider the case when the number of variables by far exceeds the number of true covariates. Therefore the cardinality of the set $\mathcal{S}$ is 
infinite, and the challenge is to find a sparse solution in it. The starting point of our analysis will be a solution of the penalized Euclidean distance problem defined by (\ref{optim}). As before, we let $\hat\theta_{j}$ represent the angle between vectors $X_{j}$ and $Y-X\hat\beta.$ We note that the angle $\hat\theta_j$ satisfies the equation
\[
\hat\theta_j=\frac{\pi}{2}-\arcsin\left(\frac{X_j^T(Y-X\hat\beta)}{\|Y-X\hat\beta\|}\right) \; \; \; \; , \; \; \; 0 \le \hat\theta_j < \pi.
\]
whenever $\|Y-X\hat\beta\|\neq 0$.
Also, let $\hat k=\left( \frac{\|\hat\beta\|}{\|\hat\beta\|_1} \right)^{1/2}$ and 
$$\frac{1}{\sqrt[4]{p}}\leq\hat k\leq 1$$ 
provided $\hat\beta\neq 0_{\mathbb{R}^p}.$ Note that $\hat k$ is a measure of sparsity, with highest value $1$ 
when there
is a single non-zero element in $\beta$ (very sparse), and with smallest value when all elements of $\beta$ are equal 
and non-zero (very non-sparse). 
We assume that $0_{\mathbb{R}^p}$ is not a minimizer of $\|Y-X\beta\|.$ 

\begin{lemma}
If $\hat\beta(\lambda)$ is a solution of (\ref{optim}) and $\hat\beta_j(\lambda)\neq 0$, then 
\begin{equation}\label{rel_mag}
\frac{\hat\beta_j(\lambda)}{\|\hat\beta(\lambda)\|}=\hat k\left(\frac{2\cos(\hat\theta_j)}{\lambda}-\hat k\sgn(\hat\beta_j(\lambda))\right)
\end{equation}  
\end{lemma}

\begin{result}
\label{iff} We have \begin{equation}\label{lam_zero}
|\cos(\hat\theta_j)|\leq\frac{\lambda\hat k}{2} \,\,\, \text{if and only if}\,\,\, \hat\beta_j(\lambda)=0.
\end{equation}
\end{result}
\begin{result}\label{sgn}
If $\hat\beta$ is the solution of (\ref{optim}) and its $j$-th component is non-zero (i.e. $\hat\beta_{j}\neq 0$) then $\sgn(\hat\beta_{j})=\sgn(X_{j}^T(Y-X\hat\beta))=\sgn(\frac{\pi}{2}-\hat\theta_{j}).$ 
\end{result}
The following result helps demonstrate the existence of a minimizing sequence whose terms have the grouping effect property for the relative size of their components.
\begin{lemma}\label{cosine}
If $\hat\beta$ is the solution of (\ref{optim}), then $\left|\frac{\hat\beta_j}{\|\hat\beta(\lambda)\|}\right|<M\leq 1$ if and only if ${\chg |\cos(\hat\theta_j)|}\leq\frac{\lambda}{2}\left(\hat k+\displaystyle\frac{M}{\hat k}\right)$, {\chg where $M$ is a constant.}  
\end{lemma}

We will use the results of this section to inform our choice of threshold for the Penalized Euclidean Distance 
numerical implementation.


\subsection{Oracle Property}
In this section we demonstrate that our method is also able to recover sparse signals without (pre)-estimates of the noise standard deviation or any knowledge about the signal. In \cite{Belloni} this property is referred as pivotal recovery. An important aspect is that an oracle theorem also brings a solid theoretical justification for the choice of the parameter $\lambda.$

We assume that $Y=X\beta^* +\sigma \epsilon$, where $\beta^*$ is the unknown true parameter value for $\beta$,  $\sigma$ is the standard deviation of the noise and $\epsilon_i$, $i=1,...,n$, are independent and identically  distributed with a normal law $\Phi_0$ with $E_{\Phi_0}(\epsilon_i)=0$ and $E_{\Phi_0}(\epsilon^2_i)=1.$ Let $\mathfrak{J}^*=\text{supp}(\beta^*)$.  
 For any candidate solution $\hat\beta$ we can use the notation $L$ for the plain Euclidean distance $L(\hat{\beta})=\|Y-X\hat{\beta}\|$ and the newly introduced norm is denoted by $\|\beta\|_{(1,2)}$, that is, $\|\beta\|_{(1,2)}= ( \|\beta\|_1\|\beta\| )^{1/2}.$

The idea behind the following considerations is the possibility of estimating the quotient $\|X^T\epsilon\|_{\infty} / \|\epsilon\|$ in probability following  \cite{Belloni}. We can use the same general result to show that the method we propose is also capable of producing pivotal recoveries of sparse signals.

Before stating the main theorem we introduce some more notation and definitions. The solution of the Penalized Euclidean Distance objective function is denoted by $\hat{\beta}(\lambda).$ 
Let $\|u\|_X$ denote $\|Xu\|$, $p^*$ the cardinality of $\mathfrak{J}^*$, $M^*=\|\beta^*\|$, $S=\|X^T\epsilon\|_{\infty}/\|\epsilon\|$, $c>1$ and, for brevity, $\bar{c}=(c+1)/(c-1)$.  Also, we write $u^*$ for the vector of components of $u$ that correspond to the non-zero $\beta^*$ elements, i.e. with indices in $\mathfrak{J}^*$. 
Also, we write ${u^*}^c$ for the vector of components of $u$ that correspond to the zero elements of  $\beta^*$, i.e. with indices in the complement of $\mathfrak{J}^*$.  We shall initially focus on the case $n^2 > p$. 
Consider
\begin{equation}
\Delta_{\bar{c}}^* =\left\{u\in\mathbb{R}^p : u\neq\{0_{\mathbb{R}^p}\},\,\,\|{u^*}^c\|_1\leq\bar{c}\|{u^*}\|_1+\frac{c\sqrt[4]{p}}{c-1}\sqrt[4]{p^*}M^*\right\}. 
\end{equation}
\noindent Assume that 
\[
\bar{k}_{\bar{c}}^*=\min\limits_{u\in \Delta_{\bar{c}}^*}\frac{1}{\sqrt{n}}\frac{\|u\|_X}{\|u\|}
\]
and
\[
k_{\bar{c}}^* = \left( 1 - \frac{1}{c} \right) \min\limits_{u\in \Delta_{\bar{c}}^*}\frac{\sqrt{p^*}\|u\|_X}{2\|{u^*}\|_1+ {\sqrt[4]{p}}\sqrt[4]{p^*}M^*}
\]
are bounded away from $0$. We make the remark that if the first compatibility condition holds, there is a relatively simple scenario when the second condition would hold, as well. If $\bar{k}_{\bar{c}}^*$ is bounded away from $0$ on $\Delta_{\bar{c}}^*$, we have that  $\| u\|_X$ must be at least 
$O(\sqrt{n})$ on $\Delta_{\bar{c}}^*$. At the same time if $\|  u^* \|$  is at most $O(p^*)$ we therefore get
\begin{equation*} 
{k}^*_{\bar{c}}  =   \sqrt{p^*} O( \sqrt{n} ) / (O(p^*) + \sqrt[4]{p} \sqrt[4]{p^*} M^* ) =  O(\sqrt{n}) / ( O(\sqrt{p^*}) + \sqrt[4]{p}  \sqrt[4]{p^*} (M^*/\sqrt{p^*}) )\end{equation*}
 and we assume  $M^* /\sqrt{p^*}$ is bounded.
 Thus, the second compatibility condition could be easily achieved in the case when $p=n^{1+\alpha_1}$ and $p^*=n^{\alpha_2}$ with $\alpha_1,\, \alpha_2>0$ and $\alpha_1+\alpha_2\leq 1.$ 
We also present later a result with certain compatibility conditions for the case when $p>n^2.$

We refer to $k_{\bar{c}}^*$ and $\bar{k}_{\bar{c}}^*$ as restricted eigenvalues. The concept of restricted eigenvalues was introduced by  \cite{Bickel} with respect to the $\ell_1$ penalty function. 
 Our definition and usage are adapted to our own objective function. 
As stated before, our oracle theorem is based on estimation of $\frac{\|X^T\epsilon\|_{\infty}}{\|\epsilon\|}.$ Directly following from Lemma 1 of \cite{Belloni}, we have:
\begin{lemma}\label{lemlam}
Given $0<\alpha<1$ and some constant $c>1$, the choice $\lambda=\frac{c\sqrt[4]{p}}{\sqrt{n}}\Phi_0^{-1}\left(1-\frac{\alpha}{2p}\right)$ satisfies $\lambda\geq c\sqrt[4]{p}S$ with probability $1-\alpha.$ 
\end{lemma}

\noindent Now we are ready to state the main result:

\begin{theorem}\label{oracle} (Signal Recovery)
Assume that $\lambda\leq\frac{\rho \sqrt[4]{p} k_{\bar{c}}^*}{\sqrt{p*}}$ for some $0<\rho<1.$ If also $\lambda\geq c\sqrt[4]{p}S$ then
\begin{equation}
(1-\rho^2) \|u\|_X  \le  \frac{2c \sqrt{p^* \log(2p/\alpha)} L(\beta^*)}{k_{\bar{c}}^* \sqrt{n}} .
\label{I1}
\end{equation}
A direct consequence is
\begin{equation*}
\bar{k}^*_{\bar{c}}\|\hat{\beta}(\lambda)-\beta^*\| = 
O\left(\frac{ \sqrt{p^* \log(2p/\alpha)}  L(\beta^*)}{(1-\rho^2)n} \right)   ,    \label{I2}
\end{equation*}
and hence if $L(\beta^*) = O_p(\sqrt{n})$ and $n \to \infty$ such that $\sqrt{p^* \log(2p/\alpha)}/\sqrt{n} \to 0$, then 
$\hat\beta(\lambda) \to \beta^*$ in probability. 
\end{theorem}

We can use the value of $\lambda$ in Lemma \ref{lemlam} for practical implementation in order to ensure $\lambda\geq c\sqrt[4]{p}S$ holds with probability $1-\alpha$. Note that the rate of convergence is asymptotically the same as rates seen in other sparse regression problems (e.g. see Negahban et al., 2012), although
as for the square root Lasso of Belloni et al. (2011) knowledge of $\sigma$ is not needed.
Also there are some circumstances when we can consider other values of $\lambda$: 

\begin{corollary} \label{oraclecor}
Let $0<\xi<1$ and 
\[
\Delta_{\xi}=\left\{u\in\mathbb{R}^p, \frac{\sqrt{n}}{\sqrt[4]{p}}\|{u^*}^c\|_1\xi\leq\|{u^*}\|_1\left(\frac{2\sqrt{n}}{\sqrt[4]{p}}-\xi\right)+\| \beta^* \|_1\left(1-\frac{\sqrt{n}}{\sqrt[4]{p}}\right)\right\}.\] 
If $k_{\xi}^* = \min\limits_{u\in \Delta_{\xi}}\frac{\frac{\sqrt{p^*}}{\sqrt{n}}\|u\|_X}{\frac{2\sqrt{n}}{\xi\sqrt[4]{p}}\|{u^*}\|_1+\frac{\| \beta^* \|_1}{\xi}\left(1-\frac{\sqrt{n}}{\sqrt[4]{p}}\right)}>k>0$ and for $\lambda=c\Phi_0^{-1}\left (1-\frac{\alpha}{2p}\right )\frac{\sqrt[4]{p}}{n}$, with $c>1$, if $\sqrt{\frac{\|\hat{\beta}(\lambda)\|}{\|\hat{\beta}(\lambda)\|_1}}-\frac{\sqrt{n}}{c\sqrt[4]{p}}\geq\xi>0$ and, at the same time, we assume $\lambda\leq\frac{\rho \sqrt[4]{p} k_{\xi}^*}{\sqrt{n}\sqrt{p^*}}$ for some $0<\rho<1$ then we also have an oracle property, i.e.
\[
\|\hat{\beta}(\lambda)-\beta^*\| =  O\left( \sqrt{ \frac{p^* \log(2p/\alpha) }{n} }  \right)
\]
with probability $1-\alpha$.
\end{corollary}
We use the corollary to suggest a method for choosing the model parameters by maximizing 
\begin{equation}
{\hat k} = \left( \frac{\|\hat{\beta}(\lambda)\|}{\|\hat{\beta}(\lambda)\|_1} \right)^{1/2} , \label{hatxi}  
\end{equation} 
which would encourage sparse models.

Although the compatibility conditions require $p \le n^2$ if we assume that $\min\limits_{u\in \Delta_{\xi}} \|u\|_X/(\sqrt{n} \|u\|_1)$ is bounded away from 0, we could allow $p > n^2$ and a sufficient condition for the compatibility inequality invoked by the Corollary would be to have $\hat k > \xi > 0$ for all $n$. Such a condition is not unrealistic if the set of columns of the design matrix $X$ has a finite partition by subsets of highly correlated covariates. If we have that the set of all indices has always a finite partition by subsets of indices whose corresponding columns in the design matrix $X$ maintain strong correlations. Thus we assume
\[  
\{1,\ldots,p\} = \displaystyle\bigcup\limits_{k=1}^{m}\mathfrak{J}_k
\]
where $m<\infty$ when $n\rightarrow\infty$. If for each subset $\mathfrak{J}_k$ the columns of $X$ whose indices are in $\mathfrak{J}_k$ have the lowest pairwise correlation, $\rho_k$, such that
\[
1-\rho_k = O\left(\frac{\sqrt{p}}{n^2|\mathfrak{J}_k|^2}\right)
\] 
then, based on Theorem \ref{grouping}, we can show that $\frac{\|\hat{\beta}(\lambda)\|_1}{\|\hat{\beta}(\lambda)\|}$ is bounded when $n\rightarrow\infty$ and implicitly the compatibility condition from the Corollary is satisfied.

For a practical implementation of our method we make use of the proven theoretical results. From the signal recovery theorem and corollary we obtain that $\|\hat{\beta}(\lambda)-\beta^*\| = O\left( \sqrt{p^* \log(2p/\alpha)}/\sqrt{n}\right)$ with probability $1-\alpha$. 
Thus, if $j$ is an index where there is no signal, i.e. $\beta_j^*=0$ then, from the previous equation, we have that $|\hat{\beta_j}(\lambda))|<\|\hat{\beta}(\lambda)-\beta^*\|\leq const. \sqrt{p^* \log(2p/\alpha)}/\sqrt{n}$. If $\|\hat{\beta}(\lambda)\|\neq 0$ we can divide by $\|\hat{\beta}(\lambda)\|$ and get 
\begin{equation}\label{deln}
\frac{|\hat{\beta_j}(\lambda))|}{\|\hat{\beta}(\lambda))\|} < \delta(p)/\sqrt{n}, 
\end{equation}
where 
\begin{equation}
\delta(p) \propto \sqrt{ {p^* \log(2p/\alpha)} } . \label{defC}
\end{equation}
We will use (\ref{deln}) to inform a threshold choice as part of the Penalized Euclidean Distance 
numerical  implementation in the next Section. As well as 
dependence on $n$ we also investigate the effect of $p$ on the relative size of the components. 
Note that the components of $\hat\beta(\lambda)$ whose relative size 
(i.e. $\frac{|\hat\beta_j(\lambda)|}{\|\hat\beta(\lambda)\|}$) is small belong to columns of $X$ that make 
with $Y-X\hat\beta(\lambda)$ an angle a lot closer to $\frac{\pi}{2}$ than the rest of the columns. For example, since 
${\chg 1/\sqrt[4]{p}}\leq\hat k\leq 1$ and if $\frac{|\hat\beta_j(\lambda)|}{\|\hat\beta(\lambda)\|}< {\chg \frac{\delta(p)}{\sqrt{n}} }$ for some $\delta(p) >0$, from equation (\ref{rel_mag}) we have
\begin{equation}\label{delta}
\cos(\hat\theta_j)<\frac{\lambda}{2}\left(\hat k+ {\chg \frac{\delta(p)}{\hat k \sqrt{n}} }\right) \le \frac{\lambda}{2}\left(\hat k+ {\chg \frac{\delta(p)p^{1/4}}{n^{1/2}} } \right),
\end{equation}
where $\lambda = O(\sqrt[4]{p}/\sqrt{n})$.

For a practical method we implement the detection of a set $I(\lambda,\delta)$ of irrelevant indices (with zero 
parameter estimates):
\begin{equation}\label{delp}
I(\lambda,\delta(p))=\bigg\{j\,\,,\,\,\frac{|\hat\beta_j(\lambda)|}{\|\hat\beta(\lambda)\|}  < {\chg \frac{\delta(p)}{\sqrt{n}} } \bigg\} ,
\end{equation}
where $\delta(p)$ is a threshold value defined in (\ref{defC}) that needs to be chosen. 
We construct a new vector $\hat{\hat{\beta}}(\lambda)$ which satisfies $\hat{\hat{\beta}}_j(\lambda)=0$, if $j\in I(\lambda,\delta(p))$ and the rest of the components give a minimizer of
\[
\|Y-\hat{\hat X}\beta\|+\lambda\left( \|\beta\|  \|\beta\|_1 \right)^{1/2} \]
where $\hat{\hat X}$ is obtained from $X$ by dropping the columns with indices in $I(\lambda,\delta(p)).$ 

In the next section  we show how these results can be used in a numerical implementation for finding sparse minimizers of $L(\beta)$.

\section{Numerical Implementation}\label{Algorithm}
The objective function $L_{PED}(\lambda,\beta)=\|Y-X\beta\|+\lambda (\|\beta\|  \|\beta\|_1)^{1/2}$ is convex for any choice of $\lambda$ and also 
differentiable on all open orthants in ${\mathbb R}^p$ bounded away from the
hyperplane $Y-X\beta = 0$. In order to find good approximations for minimizers of our objective function, as in many cases of nonlinear large scale convex optimization problems, a Quasi-Newton method may be preferred since it is known to be considerably faster than methods like coordinate descent by achieving super-linear convergence rates. Another important advantage is that second-derivatives are not necessarily required. For testing purposes, we present a 
numerical implementation based on the well performing Quasi-Newton methods for convex optimization known as  
Broyden-Fletcher-Goldfarb-Shanno (BFGS) methods: limited-memory BFGS
(L-BFGS) \citep{Nocedal:80} and BFGS \citep{Bonnans:06}. We also tested a version of non-smooth
BFGS called Hybrid Algorithm for Non-Smooth Optimization (HANSO) \citep{Lewis:08} and obtained very similar results. 

The idea for the estimation is to use theoretically informed parameters based on the Theorem 2, Corollary 1 and (\ref{delta}), in order to choose a suitable value of $\lambda$ and give a sparse estimate of $\beta^*$ after thresholding. We are choosing a lambda value in the interval between $\Phi_0^{-1}\left (1-\frac{\alpha}{2p}\right )\frac{\sqrt[4]{p}}{n}$ and $\Phi_0^{-1}\left (1-\frac{\alpha}{2p}\right )\frac{\sqrt[4]{p}}{\sqrt{n}}$ that is maximising $\hat k(\lambda_0) := \left(\frac{\|\hat{\beta}(\lambda_0)\|}{\|\hat{\beta}(\lambda_0)\|_1}\right)^{1/2}$. We retain the components of the solution that have higher relative contributions, i.e. $\frac{|\hat\beta_j|}{\|\hat\beta\|}\geq \delta(p) /\sqrt{n}$ where $\delta(p)$ is a tuning thresholding constant that could be selected by some information criterion such as AIC or by n-fold cross validations or we could fix $\delta$, e.g. $\delta = 0.75$.   
The steps for the numerical approximation of $\beta^*$ by using the Penalized Euclidean Distance method are as follows:
{
\begin{enumerate}
 \setlength{\itemsep}{1pt}
  \setlength{\parskip}{0pt}
  \setlength{\parsep}{0pt}
\item Use a Quasi-Newton algorithm (e.g. L-BFGS) to minimize the convex objective function (\ref{obj}) with $\lambda $ values between $\Phi_0^{-1}\left (1-\frac{\alpha}{2p}\right )\frac{\sqrt[4]{p}}{n}$ and $\Phi_0^{-1}\left (1-\frac{\alpha}{2p}\right )\frac{\sqrt[4]{p}}{\sqrt{n}}$ and evaluate $\hat k(\lambda)$.
\item For the solution $\hat\beta$ that maximises $\hat k(\lambda)$ ,  set $\hat\beta_j=0$ if $\frac{|\hat\beta_j|}{\|\hat\beta\|}\leq \delta(p) /\sqrt{n}$  (the choice of $\delta(p)$ is motivated by (\ref{delta})). Eliminate the columns of the design matrix corresponding to the zero coefficients $\hat\beta_j$, with $p^*$ columns remaining.
\item Use the Quasi-Newton algorithm to minimise the objective function with the remaining columns of the design matrix and $\lambda$ between $\Phi_0^{-1}\left (1-\frac{\alpha}{2p}\right )\frac{\sqrt[4]{p}}{n}$ and $\frac{\sqrt[4]{p^*}}{\sqrt{n}}$  and output the solution. 
\end{enumerate}
}

\noindent For all the numerical simulations and almost all real data sets a default value of $\lambda$ was used for the last step of the numerical approximation, namely $\frac{\sqrt[4]{p^*}}{\sqrt{n}}$.

\section{Numerical Applications}
\subsection{Simulation study}\label{Simulations}
\noindent {\it Example 1.} 
We consider a simulation study to illustrate the performance of our method and the grouping effect in the case when $p \ge n$. In this example we compare the results with the Square Root of Lasso method \citep{Belloni} that uses a scaled Euclidean distance as a loss function plus an $\ell_1$ penalty term, using the asymptotic choice of $\lambda$. We also compare the results with both Lasso and Elastic Net methods as they are implemented in the publicly available packages for R, again using the default options. 
In particular we used 10-fold cross-validation 
to choose the roughness penalty for Lasso and the Elastic Net using the command {\tt cv.glmnet} in the R package {\tt glmnet} and we use the command {\tt slim} in the R package {\tt flare} 
with penalty term 
$\lambda = 1.1 \Phi_0^{-1}(1-0.05/(2p))/\sqrt{n}$. We use the Penalized Euclidean Distance method with a default $\delta = 0.75$  or chosen with the AIC criterion from a range of values between $0.75$ and $1.5$ (\ref{hatxi}). 

We consider situations with weak signal, strong noise and various correlated designs. In particular, for a range of 
values of $n, p, \rho$ the data are generated from the linear model
$Y = X \beta^* + \sigma\epsilon,$ where 
\[
\beta^*=(\underbrace{0.3,...,0.3}_\text{4},\underbrace{0,...,0}_\text{50},\underbrace{0.3,...,0.3}_\text{4},\underbrace{0,...,0}_\text{50},\underbrace{0.3,...,0.3}_\text{4},0,...,0),
\]
$\|\beta^*\|_0=12$, $\sigma=1.5$ and $X$ generated from a $p$-dimensional multivariate normal distribution with mean 
zero and correlation matrix $\Sigma$ where the $(j,k)$-th entry of $\Sigma$ is $\rho^{|j-k|},\,\, 1\leq j,k \leq p$.

\noindent The results are summarized in the following table and the reported values are based on averaging over 100 data sets. The distance between the true signal and the solution produced is also recorded. 
Under highly correlated designs, the method we propose shows a very efficient performance against the ``curse of dimensionality" and overcrowding.

\begin{table}[htbp]
\centering
\begin{small}
\begin{tabular}{l|rrr|rrr|rrr}
              &      &    $\rho=0.5$  &  &    &  $\rho=0.9$  &      &    &  $\rho=0.99$ &       \\
& TP & MS & RMSE & TP & MS   & RMSE & TP & MS   & RMSE  \\     
\hline
$n=100,p=200$  &&& &&& &&    \\               
PED            &  {\bf 9.96}  &  20.80  &  {\bf 0.775} &  {\bf 11.62}  &  26.97 &  {\bf 0.664} & {\bf 10.17}  &  39.28  &  {\bf 0.862}\\ 
PED(AIC)       &  9.37  &  17.37  &  0.789  &  10.57  &  17.87  &  0.692 & 7.07  &  15.87  &  1.013\\ 
Elastic Net    &  9.34  &  30.17  &  0.949 &  9.43  &  25.63  &  1.064 &  6.2  &  23.52  &  1.436\\
Lasso          &  8.75  &  26.18  &   0.879 &  7.38  &  20.15  &  1.069 &  3.59  &  14.34  &  1.637\\
Sq.Rt.Lasso    &  1.18  &  1.2  &  1.015   &  3.94  &  4.65  &  0.927 &  3.03  &  8.84  &  1.291\\
\hline
$n=100,p=1000$  &&& &&& && \\
PED      &  {\bf 9.53}  &  47.34  &  1.049 & {\bf 11.52}  &  35.44  &  {\bf 0.706} &{\bf 11.14}  &  72.19  & {\bf  0.847}\\
PED(AIC) & 9.34  &  43.72  &  1.051  & 10.48  &  23.93  &  0.741 & 8.87  &  58.94  &  1.012\\
Elastic Net   & 7.93  &  44.48  &  1.042& 9.46  &  37.73  &  1.077 & 7.51  &  34.46  &  1.265 \\
Lasso   & 7.07  &  32.27  & {\bf  0.956}     & 6.87  &  27  &  1.044 & 3.3  &  18.03  &  1.510\\
Sq.Rt.Lasso    &0.9  &  0.94  &  1.025 & 2.91  &  3.42  &  0.965 & 2.97  &  9.07  &  1.253\\
\hline
$n=200,p=200$  &&& &&& && \\
PED  & {\bf 11.71}  &  21.19  &  {\bf 0.581} & {\bf 11.97}  &  27.15  &  0.540 & {\bf 11.46}  &  44.47  &  0.781\\
PED(AIC) &  11.62  &  20.31  &  0.591 &11.83  &  20.3  &  {\bf 0.535} & 10.27  &  25.09  &  {\bf 0.777}\\
Elastic Net &  11.41  &  34.3  &  0.687 & 10.45  &  25.57  &  0.860 &7.72  &  23.53  &  1.117\\
Lasso   &  11.03  &  30.89  &  0.653 &8.85  &  21.59  &  0.900 & 5.04  &  15.3  &  1.396\\
Sq.Rt.Lasso &  5.42  &  5.47  &  0.900 &7.26  &  8.3  &  0.830 &4.84  &  11.32  &  1.190\\
\hline
$n=200,p=2000$  &&& &&& && \\
PED & {\bf 11.21}  &  65.23  &  0.880 & {\bf 11.76}  &  31.79  &  {\bf 0.556} & {\bf 11.99}  &  117.42  &  {\bf 0.838}\\
PED(AIC) &  11.20  &  61.35  &  0.877 &  11.36  &  23.95  &  0.577 & 9.52  &  52.21  &  0.908  \\
Elastic Net & 10.52  &  57.22  &  0.794 &10.79  &  42.47  &  0.835  & 9.55  &  37.03  &  1.007  \\
Lasso  & 9.79  &  44.42  & {\bf  0.752} &8.65  &  31.19  &  0.856  & 4.77  &  19.74  &  1.280 \\
Sq.Rt.Lasso &  3.57  &  3.57  &  0.961 &6.45  &  7.09  &  0.837 & 4.8  &  11.14  &  1.155 \\
\hline
$n=200,p=3000$  &&& &&& && \\
PED & {\bf 11.17}  &  73.91  &  0.977 & {\bf 11.79}  &  35.04  &  {\bf 0.581} & {\bf 11.7}  &  121.92  &  {\bf 0.854}\\
PED(AIC) &  11.16  &  72.10  &  0.975 &  11.66  &  27.69  &  0.586 & 9.39  &  77.3.02  &  0.934\\
Elastic Net & 10.23  &  64.83  &  0.834 & 10.75  &  44.16  &  0.839 & 9.94  &  37.42  &  0.993\\
Lasso & 9.61  &  49.97  &  {\bf 0.788} & 8.48  &  33.77  &  0.870 &  4.62  &  19.04  &  1.293\\
Sq.Rt.Lasso & 2.86  &  2.87  &  0.980 & 6.31  &  7.05  &  0.843 &  4.6  &  10.79  &  1.173\\
\hline 
\end{tabular}
\vskip 0.75cm
\end{small}
\caption{Simulation results based on Example 1. 
The True Positives (TP) are the average number of non-zero parameters which are estimated as non-zero and the Model Size (MS) is the average number estimated non-zero parameters, 
 from 100 simulations. The Root Mean Square Error (RMSE) is given for estimating $\beta^*$. The best values in the TP and RMSE columns are in bold. }
\label{table:simu1}
\end{table}

Penalized Euclidean Distance (indicated by PED in the table) has performed very well obtaining the highest rate of true positives in many examples. 
Also, we compare with Penalized Euclidean Distance when the parameters are selected
by Akaike's Information Criterion, given by PED(AIC) in the table which also performs well. The Elastic Net is the next best, and 
Lasso (for the strongly correlated case) and Square root Lasso have low rates of True Positives. Note that the Square Root Lasso has performed rather differently here from the others. It is the only method using an asymptotic value of $\lambda$, where $n$ may not be large enough here.    
Penalized Euclidean Distance has a lower model size compared to the Elastic Net. Finally the Root Mean Square Error is generally best for Penalized Euclidean Distance, particularly for
the higher correlated situations. Overall Penalized Euclidean Distance has performed extremely well in these simulations.  

The next table summarises the True Positive rates for fixed model sizes. The fixed size model means retaining the prescribed number of top largest relative contributions:

\begin{table}[htbp]
\centering
\begin{small}
\begin{tabular}{l|rrr|rrr|rrr}
              &      &    $\rho=0.5$  &  &    &  $\rho=0.9$  &      &    &  $\rho=0.99$ &       \\
& MS=10 & MS=20 & MS=30 & MS=10 & MS=20   & MS=30 & MS=10 & MS=20   & MS=30  \\     
\hline
$n=100,p=200$  &&& &&& &&    \\               
PED            &  {\bf 7.10}  &  {\bf 9.86}  &  {\bf 11.21} &  {\bf 6.79}  & {\bf 9.81}&  {\bf 12} & {\bf 4.98}  &  {\bf 7.33}  &  {\bf 8.26}\\ 
Elastic Net    &  6.54  &  8.55  &  9.53 &  6.51  &  8.27  &  9.04 &  3.02  &  5  &  6.4\\
Lasso          &  6.13  &  8.14  &  9.16 &  5.39  &  6.75  &  7.39 &  2.81  &  3.42  &  4.05\\
Sq.Rt.Lasso    &  6.17  &  8.18  &  9.04   &  6.69  &  7.52  &  7.58 &  3.33  &  4.02  &  4.05\\
{Sq.Rt.Elastic Net} & 3.78 & 5.73 & 6.06 & 4.33 & 6.31 & 7.68 & 3.85 & 5.39 & 6.38 \\
\hline
$n=100,p=1000$  &&& &&& && \\
PED      &  {\bf 6.88}  &  {\bf 8.94}  &  {\bf 10.24} & {\bf 6.84}  &  {\bf 9.39}  &  {\bf 11.34} &{\bf 4.61}  &  {\bf 7.18}  & {\bf  8.15}\\
Elastic Net   & 5.37  &  6.98  &  7.84 & 5.37  &  6.18  &  6.93 & 3.18  &  5.12  &  6.51 \\
Lasso   & 5.11  &  6.56  &   7.20  & 5.42  &  6.27  &  6.66 & 2.85  &  3.43  &  3.75\\
Sq.Rt.Lasso    & 5.75  &  7.01  &  7.43 & 5.37  &  6.18  &  6.76 & 3.39  &  3.67  &  3.9\\
{Sq.Rt.Elastic Net} & 1.09 & 1.69 & 2.17 & 2.12 & 2.99 & 3.59 & 2.10 & 2.95 & 3.67 \\
\hline
$n=200,p=200$  &&& &&& && \\
PED  & {\bf 8.82}  &  {\bf 11.11}  &  {\bf 11.59} & 7.51  &  {\bf 10.27}  &  {\bf 11.95} & {\bf 6.57}  & {\bf 9.62}  &  {\bf 11.17}\\
Elastic Net &  8.19  &  9.95  &  10.39 & {\bf 7.82}  &  9.59  &  10.19 & 4.88  &  7.09  &  8.5\\
Lasso   &  7.64  &  9.21  &  9.62 & 7.14  &  8.27  &  8.93 & 3.5  &  4.37  &  4.96\\
Sq.Rt.Lasso &  8.53  &  10.55  &  11.11 &7.00  &  8.32  &  8.85 & 3.30  &  4.41  &  5.19\\
{Sq.Rt.Elastic Net} & 5.96 & 7.54 & 8.55 & 6.32 & 8.68 & 10.03 & 5.08 & 7.23 & 8.75\\
\hline
$n=200,p=2000$  &&& &&& && \\
PED & {\bf 7.92}  &  {\bf 9.60}  &  {\bf 10.04} & 7.47  &  {\bf 10.19}  &  {\bf 12} & {\bf 6.08}  &  {\bf 8.76}  &  {\bf 9.74}\\
Elastic Net & 7.03  &  8.34  &  8.95 & {\bf 7.85}  &  9.48  &  9.94  & 4.08  &  5.98  &  7.41  \\
Lasso  & 6.96  &  8.45  & 8.99 & 6.78  &  7.78  &  8.15  & 3.91  &  4.29  &  4.57 \\
Sq.Rt.Lasso &  7.47  &  9.01  &  9.72 & 7.02  &  8.25  &  8.77 & 3.89  &  4.52  &  4.72 \\
{Sq.Rt.Elastic Net} & 1.54 & 2.16 & 2.61 & 2.27 & 3.57 & 4.51 & 1.97 & 2.91 & 3.80\\
\hline
$n=200,p=3000$  &&& &&& && \\
PED & {\bf 8.99}  & {\bf 11.35} &  {\bf 11.95} & 7.2  &  {\bf 10.5}   &  {\bf 12} & {\bf 6.77}  & {\bf 9.24} &  {\bf 10.07}\\
Elastic Net & 7.36  &  8.6  &  9.09 & {\bf 8.13}  &  9.89  &  10.21 & 4.46  &  6.65  &  8.14\\
Lasso & 7.13  &  8.32  &  8.66 & 6.49  &  7.31  &  7.61 &  4.09  &  4.61  &  4.84\\
Sq.Rt.Lasso & 7.35  &  8.86  &  9.4 & 7.11  &  8.15  &  8.51 &  4.08  &  4.65  &  4.89\\
{Sq.Rt.Elastic Net} & 0.46 & 0.79 & 1.21 & 1.05 & 1.92 & 2.68 & 1.46 & 2.22 & 2.82\\
\hline 
\end{tabular}
\vskip 0.75cm
\end{small}
\caption{Simulation results based on Example 1 and the numerical implementation described in the previous section. The best values for the True Positives  are in bold. 
The True Positives (TP) are the average number of non-zero parameters which are estimated as non-zero and the fixed Model Size (MS) 
 from 100 simulations. The "Sq. Rt. Elastic Net" refers to minimizing the square root loss function with the Elastic Net penalty with $\lambda =  \Phi_0^{-1}(1-0.05/(2p))\frac{p^{1/4}}{\sqrt{n}}.$}
\label{table:simu2}
\end{table}

\subsection{Real Data Applications}\label{Numerical Applications}

\noindent Description of the data sets:

\vspace{3mm}

\noindent \textit{Air:} Daily air quality measurements in New York, May to September 1973. Data consists of 5 variables with 111 observations after removing missing values (see, Chambers \textit{et al.}, 1983).\\

\noindent \textit{Servo:} According to Qunilan (1993), this data was collected from a simulation of a servo system involving a servo amplifier, a motor, a lead screw/nut, and some sort of sliding carriage (see, Quanlan, 1993). The data consists of 167 observations and 4 factors. Factors were converted to dummy variables, for a total of 19 variables. The data is available at the UCI repository (see, Asuncion and Newman, 2007).\\

\noindent \textit{Tecator:} The goal for collecting this data set was to predict the fat content of a meat sample on the basis of its near infrared absorbance spectrum. The Tecator2 data set consists of 215 observations and 100 variables comprising log absorption values. The 215 observations are the training (C), monitoring (M), and testing (T) data described in Borggaard and Thodberg (1992).\\

\noindent \textit{Housing:} Median house price for 506 census tracts of Boston from the 1970 census. The data comprises 506 observations and 13 variables. For details, see Harrison and Rubinfeld (1978). The data is available at the UCI repository (see, Asuncion and Newman, 2007).\\

\noindent \textit{Ozone:} Los Angeles ozone pollution data from 1976. After removing missing values, the data set comprised 203 observations on 12 variables, each observation is one day. For details, see Breiman and Friedman (1985).\\

\noindent \textit{Iowa:} The Iowa wheat yield data consists of 33 observations and 9 variables. For details, see CAED report (1998).

It is important to note that throughout the simulations and the real data analyses both Lasso and Elastic Net were run with double cross validation for selecting the model size and the tuning parameter. We used Penalized Euclidean Distance with default values and also we ran Penalized Euclidean Distance with a single 10-fold cross validation for tuning $\delta$ that affects only the model size; the results are reported in the table below. 
In some cases, such as the Tecator data set, the prediction error further improved when the final choice of  $\lambda$ was $\Phi_0^{-1}\left (1-\frac{\alpha}{2p}\right )\frac{\sqrt[4]{p}}{n}.$ Also, in the case of the Servo data set, the variable selection benefitted from searching for lambda in a subinterval of the one proposed by default, namely between $\Phi_0^{-1}\left (1-\frac{\alpha}{2p}\right )\frac{\sqrt[4]{p}}{n}$ and $\frac{\sqrt[4]{p}}{\sqrt{n}}.$


\begin{table}[htbp]
\centering
\caption{The table presents the mean square error rates in 90 - 10 cross validations}
\label{Real Data Applications}
\begin{small}
\begin{tabular}{|
>{\columncolor[HTML]{EFEFEF}}l |r|r|r|r|r|}
\hline
         & \multicolumn{1}{l|}{\cellcolor[HTML]{EFEFEF}PED} & \multicolumn{1}{l|}{\cellcolor[HTML]{EFEFEF}PED-CV} & \multicolumn{1}{l|}{\cellcolor[HTML]{EFEFEF}Elastic Net} & \multicolumn{1}{l|}{\cellcolor[HTML]{EFEFEF}Lasso} & \multicolumn{1}{l|}{\cellcolor[HTML]{EFEFEF}Sqrt. Lasso} \\ \hline
Servo  & 0.3436 &    0.2977                                           & 0.3041                                                   	& 0.3130                                             &        0.7393                                                  \\ \hline
Pollute &2334.0 &  2197.3                                           & 2334.8                                                     & 2778.6                                             &        3066.9                                                  \\ \hline
Iowa   &89.3389 & 100.8377                                       & 143.8566                                                & 122.8793                                           &     133.9454                                                     \\ \hline
Air     &0.2348 & 0.2413                                             & 0.2811                                                       & 0.2874                                             &           0.3962                                               \\ \hline
Ozone  &17.8964 & 15.6248                                       & 16.8717                                                   & 17.1563                                            &     21.4246                                                     \\ \hline
Tecator &12.0489 & 10.8289                                          & 50.0329                                                 & 44.3648                                            &     131.1065                                                     \\ \hline
\end{tabular}
\end{small}
\end{table}

\noindent \textit{Melanoma:} In this application we implement Penalized Euclidean Distance as a variable selection tool when the response variable serves for binary classification. 
We consider an application of the method to a proteomics dataset from the study of melanoma (skin cancer). The mass spectrometry dataset was described by \citet{Mian et al 2005} and 
further analysed by \citet{Browne:2010}. The data consist of mass spectrometry scans  from serum samples of $205$ patients, 
with $101$ patients with Stage I melanoma (least severe) and $104$ patients with Stage IV melanoma (most severe). Each mass spectrometry scan 
consists of an intensity for $13951$ mass over charge ($m/z$) values between $2000$ and $30000$ Daltons. It is of interest to
find which m/z values could be associated with the stage of the disease, which could point to potential proteins for use as biomarkers.
We first fit a set of 500 important peaks to the overall mean of the scans using the deterministic peak finding algorithm of \citet{Browne:2010} to 
obtain 500 m/z values at peak locations.  
We consider the disease stage to be the response, with $Y=-1$ for Stage I and $Y=1$ for Stage IV. Note that we have an ordered 
response here as Stage IV is much more severe than Stage I, and it is reasonable to treat the problem as a regression problem.

We fit the Penalized Euclidean Distance regression model versus the 
intensities at the 500 peak locations. We have $n=205$ by $p=500$.
The data are available at {\tt http://www.maths.nottingham.ac.uk/$\sim$ild/mass-spec}

Here we use $\alpha=0.05$. The parameter values chosen 
to maximize $\hat k $ are $\lambda=0.5$ and $\delta(p)=0.75$, selecting 96 non-zero m/z values. 
\citet{Browne:2010} also considered 
a mixed effects Gaussian mixture model and a two stage t-test for detecting significant peaks. 
If we restrict ourselves to the coefficients corresponding to the 50 largest peaks in height, \citet{Browne:2010} identified 
17 as non-zero as did Penalized Euclidean Distance, with 8 out of the 17 in common. 
If we apply PED(AIC) then 7 peaks are chosen out of the largest 50 of which only 2 are in common 
with \citet{Browne:2010}. The Elastic Net chose 6 peaks with 5 of those in common with \citet{Browne:2010} 
and for the Lasso 5 peaks were chosen from the top 50 largest, with 4 in common with \citet{Browne:2010}. 
Note that here Penalized Euclidean Distance has selected the most peaks in common with \citet{Browne:2010}, and it is 
reassuring that the different methods have selected some common peaks.

\section*{Acknowledgement}
We would like to thank the editor, associate editor and the referees for their valuable comments and insight that significantly helped to improve the quality of this manuscript.


\vskip 2cm

\section*{Appendix A: Proofs}
{\begin{proof} (Lemma \ref{convlema}.)
\noindent Let $${\cal C}=\{\beta\in\mathbb{R}^n : \sqrt{f_{p_1}(\beta)  f_{p_2}(\beta)}\leq 1\}.$$  We note that
$
{\cal C}\equiv\{\beta\in\mathbb{R}^n : [\sqrt{f_{p_1}(\beta)  f_{p_2}(\beta)}]^{2p_1p_2}\leq 1\}
$
and therefore ${\cal C}$ is a bounded, closed and convex subset of $\mathbb{R}^n$ which contains the origin. Let $g(\beta)=\left( \displaystyle\sum\limits_{i=1}^{n}|\beta_i|^{p_1}\right)^{p_2}\left( \displaystyle\sum\limits_{i=1}^{n}|\beta_i|^{p_2}\right)^{p_1}$ and let $\Epi(g)$ denote its epigraph, i.e. $\Epi(g)=\{ (\beta,t)\in\mathbb{R}^{n+1} \, : \, g(\beta)\leq t \}.$ The set  ${\cal C}$ is convex and orthant symmetric.  Indeed, the Hessian of $g$ is positive semi-definite on each orthant of $\mathbb{R}^n$ since, after differentiating $g$ twice with the product rule, it can be written as a sum of three matrices which can be argued, by applying Sylvester's Theorem, that is positive semi-definite. 

We see that in our case $\Epi(f_{p_1}  f_{p_2})=\{t({\cal C},1) : \, t\in [0,+\infty)\}$ and therefore $\Epi(f_{p_1}  f_{p_2})$ is a convex cone in $\mathbb{R}^{n+1}$ since ${\cal C}$ is a convex set  in $\mathbb{R}^n$. This shows that $\sqrt{f_{p_1}  f_{p_2}}$ is a convex function. Because $\sqrt{f_{p_1}  f_{p_2}}$ is convex and homogeneous of degree $1$ it follows that it must also satisfy the triangle inequality. Therefore $\sqrt{f_{p_1}  f_{p_2}}$ is a norm on $\mathbb{R}^n.$ 
\end{proof}


\begin{proof} (Theorem \ref{grouping}.) 
Since $\hat\beta(\lambda)=\arg\min\limits_{\beta}\{L_{PED}(\lambda,\beta)\}$ we have 
\begin{equation}
\frac{\partial L_{PED}(\lambda,\beta)}{\partial \beta_k}\bigg|_{\beta=\hat\beta(\lambda)}=0\,\,\,\text{for every}\,k=1,2,...p
\end{equation}
unless $\hat\beta_k(\lambda)=0.$ Thus, if $\hat\beta_k(\lambda)\neq 0$ we have
\begin{equation}
-\frac{X_k^T[Y-X\hat\beta(\lambda)]}{\|Y-X\hat\beta(\lambda)\|}+\frac{\lambda}{2}\frac{\frac{\hat\beta_k(\lambda)}{\|\hat\beta(\lambda)\|}|\hat\beta(\lambda)|_1}{\sqrt{\|\hat\beta(\lambda)\|  |\hat\beta(\lambda)|_1}}+\frac{\lambda}{2}\frac{{\chg \sgn\{\hat\beta_k(\lambda)\} }\|\hat\beta(\lambda)\|}{\sqrt{\|\hat\beta(\lambda)\|  |\hat\beta(\lambda)|_1}}=0.
\end{equation}
If we take $k=i$ and $k=j$, after subtraction we obtain
\begin{equation}
\frac{[X_j^T - X_i^T][Y-X\hat\beta(\lambda)]}{\|Y-X\hat\beta(\lambda)\|}+\frac{\lambda}{2}\frac{[\hat\beta_i(\lambda)-\hat\beta_j(\lambda)]|\hat\beta(\lambda)|_1}{\sqrt{\|\hat\beta(\lambda)\|^3  |\hat\beta(\lambda)|_1}}=0
\end{equation}
since $\sgn\{\hat\beta_i(\lambda)\}=\sgn\{\hat\beta_j(\lambda)\}.$ Thus we get
\begin{equation}\label{ineq1}
\frac{\hat\beta_i(\lambda)-\hat\beta_j(\lambda)}{\|\hat\beta(\lambda)\|}=\frac{2}{\lambda}\frac{\sqrt{\|\hat\beta(\lambda)\|  |\hat\beta(\lambda)|_1}}{|\hat\beta(\lambda)|_1}[X_j^T-X_i^T]\hat r(\lambda)
\end{equation}
where $\hat r(\lambda)=\frac{y-X\hat\beta(\lambda)}{\|y-X\hat\beta(\lambda)\|}$ and
$
{\chg \|X_j^T - X_i^T\|^2 }=2(1-\rho)
$
since $X$ is standardized, and {\chg $\rho = \cos(\theta_{ij})$}. We have
$
\frac{\sqrt{\|\beta\|  \|\beta\|_1}}{\|\beta\|_1}\leq 1
$
for any nonzero vector $\beta$ in $\mathbb{R}^p$ and
$
|\hat r(\lambda)|\leq 1.
$
Thus, equation (\ref{ineq1}) implies that
\begin{equation}
D_\lambda(i,j)\leq\frac{{\chg 2}|\hat r(\lambda)|}{\lambda}{\chg \|X_i - X_j\|}\leq\frac{{\chg 2}}{\lambda}\sqrt{2(1-\rho)} {\chg \le 2 \frac{\theta_{ij}}{\lambda},}
\end{equation}
which proves the grouping effect property for the proposed method.
\end{proof}

\begin{proof}(Proposition \ref{iff})
Here we are going to prove the necessity part of the statement since the sufficiency follows directly from the previous Lemma.
Let us assume that 
\[
\hat\beta(\lambda)=(\hat\beta_1(\lambda),...,\hat\beta_{j-1}(\lambda),0,\hat\beta_{j+1}(\lambda)...\hat\beta_p(\lambda))=\arg\min\limits_{\beta}\{L_{PED}(\lambda,\beta)\}
\]
for a given $\lambda>0.$ Here we can fix $\lambda$ and, for brevity, we can omit it from notations in the course of this proof. For any $t>0$ we have
\[
\frac{L_{PED}(\hat\beta_1,...\hat\beta_{j-1},t,\hat\beta_{j+1},...\hat\beta_p)-L_{PED}(\hat\beta_1,...\hat\beta_{j-1},0,\hat\beta_{j+1},...\hat\beta_p)}{t}\geq 0.
\]
Again, for brevity we can denote $\hat\beta_{t@j}=(\hat\beta_1,...\hat\beta_{j-1},t,\hat\beta_{j+1},...\hat\beta_p)^T$ and also let $\hat\theta_{t@j}$ be the angle between $x_{\ast,j}$ and $Y-X\hat\beta_{t@j}.$ By using  the mean value theorem (Lagrange), there exists $0<t^*<t$ such that 
\[
\frac{L_{PED}(\hat\beta_{t@j})-L_{PED}(\hat\beta_{0@j})}{t}=-\cos(\hat\theta_{t^*@j})+\lambda\frac{\sqrt{\|\hat\beta_{t@j}\|  |\hat\beta_{t@j}|_1}-\sqrt{\|\hat\beta_{0@j}\|  |\hat\beta_{0@j}|_1}}{t}
\] 
If we rationalize the numerator of the second fraction in the previous equation, we get
\[
\frac{L_{PED}(\hat\beta_{t@j})-L_{PED}(\hat\beta_{0@j})}{t}=-\cos(\hat\theta_{t^*@j})+\lambda\frac{\frac{\|\hat\beta_{t@j}\|  |\hat\beta_{t@j}|_1-\|\hat\beta_{0@j}\|  |\hat\beta_{0@j}|_1}{t}}{\sqrt{\|\hat\beta_{t@j}\|  |\hat\beta_{t@j}|_1}+\sqrt{\|\hat\beta_{0@j}\|  |\hat\beta_{0@j}|_1}}
\]
and thus
\[
\cos(\hat\theta_{t^*@j})\leq \lambda\frac{\frac{\|\hat\beta_{t@j}\|  |\hat\beta_{t@j}|_1-\|\hat\beta_{0@j}\|  |\hat\beta_{0@j}|_1}{t}}{\sqrt{\|\hat\beta_{t@j}\|  |\hat\beta_{t@j}|_1}+\sqrt{\|\hat\beta_{0@j}\|  |\hat\beta_{0@j}|_1}}.
\]
Also
\[
\frac{\|\hat\beta_{t@j}\|  |\hat\beta_{t@j}|_1-\|\hat\beta_{0@j}\|  |\hat\beta_{0@j}|_1}{t}=|\hat\beta_{t@j}|_1\frac{\|\hat\beta_{t@j}\|-\|\hat\beta_{0@j}\|}{t}+\|\hat\beta_{0@j}\|\frac{|\hat\beta_{t@j}|_1-|\hat\beta_{0@j}|_1}{t}
\]
and we notice that $\frac{|\hat\beta_{t@j}|_1-|\hat\beta_{0@j}|_1}{t}=1$ for any $t>0.$ Letting $t\rightarrow 0$ we obtain
\[
\cos(\hat\theta_{0@j})\leq \frac{\lambda}{2}\sqrt{\frac{\|\hat\beta_{0@j}\|}{|\hat\beta_{0@j}|_1}}=\frac{\lambda\hat k}{2}.
\]
Analogously, by starting with $t<0$, we can show that 
\[
\cos(\hat\theta_{0@j})\geq -\frac{\lambda}{2}\sqrt{\frac{\|\hat\beta_{0@j}\|}{|\hat\beta_{0@j}|_1}}=\frac{\lambda\hat k}{2}.
\]

\end{proof}

\begin{proof}(Proposition \ref{sgn}) By writing the necessary conditions for optimality in the case of problem (\ref{optim}) we have 
\[
\sgn(X_{j}^T(Y-X\hat\beta))=\sgn\left(\frac{\pi}{2}-\hat\theta_{j}\right)
\]
and
\[
\frac{\hat\beta_j(\lambda)}{\|\hat\beta(\lambda)\|}=\hat k\left(\frac{{\chg 2} X_{j}^T(Y-X\hat\beta)}{{\chg \lambda} \|Y-X\hat\beta\|}- {\sgn(\hat\beta_j(\lambda))}{\hat k}\right)
\]
if $\hat\beta_j(\lambda)\neq 0.$ Since $\hat k>0$ we have $\sgn(\hat\beta_{j})=\sgn(X_{j}^T(Y-X\hat\beta))=\sgn(\cos(\hat\theta_j)).$ 
\end{proof}

\begin{proof}(Lemma \ref{cosine})
The proof follows directly from (\ref{rel_mag}) and (\ref{lam_zero}).
\end{proof}

\noindent We make the observation that if $\hat\beta(\lambda)$ is a solution of (\ref{optim}) we have $\cos(\hat\theta_j)\leq\frac{\lambda}{2}\left(\hat k+\displaystyle\frac{M}{\hat k}\right)$ and therefore $\cos(\hat\theta_j)\rightarrow 0$ when $\lambda\rightarrow 0$ since $M\leq 1$ and $p^{-1/4}\leq\hat k\leq 1.$

{
{\it Proof.} {(Theorem \ref{oracle})} The proof follows a similar method to that of Theorem 1 in \cite{Belloni}.
Given that $\hat\beta(\lambda)$ is a minimizer of the PED objective function for a given $\lambda$, we have 
\[
L(\hat{\beta}(\lambda))-L(\beta^*)\leq \lambda \| \beta^* \|_1-\lambda\|\hat{\beta}(\lambda)\|_{(1,2)}\leq \lambda \| \beta^* \|_1-\frac{\lambda}{\sqrt[4]{p}}\|\hat{\beta}(\lambda)\|_1.
\]
We obtain
\[
L(\hat{\beta}(\lambda))-L(\beta^*)\leq \frac{\lambda}{\sqrt[4]{p}}\| \beta^* \|_1-\frac{\lambda}{\sqrt[4]{p}}\|\hat{\beta}(\lambda)\|_1+\sqrt[4]{p}\lambda M^* \leq \frac{\lambda}{\sqrt[4]{p}}(\|{u^*}\|_1-\|{u^*}^c\|_1)+\lambda M^*\sqrt[4]{p^*}.
\]
At the same time, due to the convexity of $L$, we have
\[
L(\hat{\beta}(\lambda))-L(\beta^*)\geq (\nabla L(\beta^*))^Tu \geq -\frac{\|X^T\epsilon\|_{\infty}}{\|\epsilon\|}\|u\|_1\geq-\frac{\lambda}{c\sqrt[4]{p}}(\|{u^*}\|_1+\|{u^*}^c\|_1)
\]
if $\lambda \ge c\sqrt[4]{p}S$, where $S=\frac{\|X^T\epsilon\|_{\infty}}{\|\epsilon\|}$.  Thus we have 
\[
\|{u^*}^c\|_1\leq\frac{c+1}{c-1}\|{u^*}\|_1+\frac{c\sqrt[4]{p}}{c-1}\sqrt[4]{p^*}M^*
\]
and also
\[
\|u\|_1\leq \frac{2c}{c-1}\|{u^*}\|_1+\frac{c\sqrt[4]{p}}{c-1}\sqrt[4]{p^*}M^*.
\]

Now 
\begin{eqnarray*} 
L(\hat{\beta}(\lambda))-L(\beta^*) & \le & | L(\hat{\beta}(\lambda))-L(\beta^*) | \le \frac{\lambda}{c\sqrt[4]{p}}(\|{u^*}\|_1+\|{u^*}^c\|_1)\\
& \le & \frac{\lambda \|{u^*}\|_1 } {c\sqrt[4]{p}} \le \frac{\lambda \sqrt{p^*} \| u \|_X }{ \sqrt[4]{p} k_{\bar{c}}^* } .
\end{eqnarray*}

Considering the identity
\[
L^2(\hat{\beta}(\lambda))-L^2(\beta^*)=\|u\|^2_X-2(\sigma\epsilon^T Xu)
\]
along with
\[
L^2(\hat{\beta}(\lambda))-L^2(\beta^*)=(L(\hat{\beta}(\lambda))-L(\beta^*))(L(\hat{\beta}(\lambda))+L(\beta^*))
\]
and the fact that
\[
2|\sigma\epsilon^T Xu|\leq 2 L(\beta^*)S\|u\|_1
\] 
we deduce
\begin{eqnarray*} 
\|u\|^2_X & \leq & \frac{\lambda\sqrt{p^*}\|u\|_X}{c\sqrt[4]{p} k_{\bar{c}}^*}\left(L(\beta^*)+\frac{\lambda\sqrt{p^*}\|u\|_X}{\sqrt[4]{p}k_{\bar{c}}^*}\right)+
  L(\beta^*)\frac{\lambda \sqrt{p^*}\|u\|_X}{\sqrt[4]{p} k_{\bar{c}}^*} ,\\
& \leq &  \frac{2\lambda\sqrt{p^*}\|u\|_X}{\sqrt[4]{p} k_{\bar{c}}^*} + \left( \frac{\lambda\sqrt{p^*}\|u\|_X}{\sqrt[4]{p}k_{\bar{c}}^*}\right)^2 . 
\end{eqnarray*}

Thus we have
\[
\left[1-\left( \frac{\lambda\sqrt{p^*}}{\sqrt[4]{p}  k_{\bar{c}}^*}\right)^2\right]\|u\|^2_X\leq \frac{2\lambda\sqrt{p^*}}{\sqrt[4]{p} k_{\bar{c}}^*}L(\beta^*)\|u\|_X  .
\]
We write $\lambda \le \frac{\rho\sqrt{p^*}}{\sqrt[4]{p}k_{\bar{c}}^*}$ where  $0 <  \rho < 1$ and 
 $\lambda=\frac{c\sqrt[4]{p}}{\sqrt{n}}\Phi_0^{-1}\left(1-\frac{\alpha}{2p}\right)$,  
and we can use the result that $\Phi_0^{-1}\left(1-\frac{\alpha}{2p}\right) \le \sqrt{ 2 \log(2p/\alpha) }$ (from Belloni et al., 2011).
Hence,  $$(1-\rho^2) \|u\|_X  \le  \frac{2c \sqrt{p^*} \sqrt{ \log(2p/\alpha)} L(\beta^*)}{k_{\bar{c}}^* \sqrt{n}}   , $$
and so
$$
\bar{k}^*_{\bar{c}}\|\hat{\beta}(\lambda)-\beta^*\|\leq 
\frac{ const. \sqrt{p^* \log(2p/\alpha)} L(\beta^*) }{(1-\rho^2){n}}   .
$$

\begin{proof}{(Corollary \ref{oraclecor})}
We have
\[
\begin{split}
L(\hat{\beta}(\lambda))-L(\beta^*)\leq \frac{\lambda\sqrt{n}}{\sqrt[4]{p}} \| \beta^* \|_1-\frac{\lambda\sqrt{n}}{\sqrt[4]{p}} \|\hat{\beta}(\lambda)\|_1+ \lambda\left( \frac{\sqrt{n}}{\sqrt[4]{p}}\|\hat{\beta}(\lambda)\|_1-\|\hat{\beta}(\lambda)\|_{(1,2)}\right)+\\
\lambda\| \beta^* \|_1\left(1-\frac{\sqrt{n}}{\sqrt[4]{p}}\right).
\end{split}
\]

\noindent If $\frac{\lambda\sqrt{n}}{c\sqrt[4]{p}}\geq S$ with probability $1-\alpha$ for some $c>1$, we obtain 
\begin{equation}\label{mainineq}
\begin{split}
-\frac{\sqrt{n}}{c\sqrt[4]{p}}\left(\|{u^*}\|_1+\|{u^*}^c\|_1\right)\leq \frac{\sqrt{n}}{\sqrt[4]{p}}\|{u^*}\|_1-\frac{\sqrt{n}}{\sqrt[4]{p}}\|{u^*}^c\|_1+\left(\frac{\sqrt{n}}{\sqrt[4]{p}}\|\hat{\beta}(\lambda)\|_1-\|\hat{\beta}(\lambda)\|_{(1,2)}\right)+\\
\| \beta^* \|_1\left(1-\frac{\sqrt{n}}{\sqrt[4]{p}}\right).
\end{split}
\end{equation}
\noindent At the same time we can write
\[
\begin{split}
\left(\frac{\sqrt{n}}{\sqrt[4]{p}}\|\hat{\beta}(\lambda)\|_1-\|\hat{\beta}(\lambda)\|_{(1,2)}\right)+\| \beta^* \|_1\left(1-\frac{\sqrt{n}}{\sqrt[4]{p}}\right)=\|\hat{\beta}(\lambda)\|_1\left(\frac{\sqrt{n}}{\sqrt[4]{p}}-\sqrt{\frac{\|\hat{\beta}(\lambda)\|}{\|\hat{\beta}(\lambda)\|_1}}\right)+\\
\| \beta^* \|_1\left(\sqrt{\frac{\|\hat{\beta}(\lambda)\|}{\|\hat{\beta}(\lambda)\|_1}}-\frac{\sqrt{n}}{\sqrt[4]{p}}\right)+\| \beta^* \|_1\left(1-\frac{\sqrt{n}}{\sqrt[4]{p}}\right)
\end{split}
\]
and thus we have 

\[
\begin{split}
\left(\frac{\sqrt{n}}{\sqrt[4]{p}}\|\hat{\beta}(\lambda)\|_1-\|\hat{\beta}(\lambda)\|_{(1,2)}\right)+\| \beta^* \|_1\left(1-\frac{\sqrt{n}}{\sqrt[4]{p}}\right)\leq\|u\|_1\left(\frac{\sqrt{n}}{\sqrt[4]{p}}-\sqrt{\frac{\|\hat{\beta}(\lambda)\|}{\|\hat{\beta}(\lambda)\|_1}}\right)+\\
\| \beta^* \|_1\left(1-\frac{\sqrt{n}}{\sqrt[4]{p}}\right).
\end{split}
\]
By combining with inequality (\ref{mainineq}) we obtain

\[
\begin{split}
\frac{\sqrt{n}}{\sqrt[4]{p}}\left(1-\frac{1}{c}\right)\|{u^*}^c\|_1\leq \frac{\sqrt{n}}{\sqrt[4]{p}}\|{u^*}\|_1\left(1+\frac{1}{c}\right)+\|u\|_1\left(\frac{\sqrt{n}}{\sqrt[4]{p}}-\sqrt{\frac{\|\hat{\beta}(\lambda)\|}{\|\hat{\beta}(\lambda)\|_1}}\right)+\\
\| \beta^* \|_1\left(1-\frac{\sqrt{n}}{\sqrt[4]{p}}\right)
\end{split}
\]
and it immediately follows that

\[
\|u\|_1\left(\sqrt{\frac{\|\hat{\beta}(\lambda)\|}{\|\hat{\beta}(\lambda)\|_1}}-\frac{\sqrt{n}}{c\sqrt[4]{p}}\right)\leq \frac{2\sqrt{n}}{\sqrt[4]{p}}\|{u^*}\|_1+\| \beta^* \|_1\left(1-\frac{\sqrt{n}}{\sqrt[4]{p}}\right).
\]

If $\sqrt{\frac{\|\hat{\beta}(\lambda)\|}{\|\hat{\beta}(\lambda)\|_1}}-\frac{\sqrt{n}}{c\sqrt[4]{p}}\geq\xi>0$ we have
\[
\|u\|_1\leq\frac{2\sqrt{n}}{\xi\sqrt[4]{p}}\|{u^*}\|_1+\frac{\| \beta^* \|_1}{\xi}\left(1-\frac{\sqrt{n}}{\sqrt[4]{p}}\right).
\]
and equivalently
\[
\|{u^*}^c\|_1\xi\leq\|{u^*}\|_1\left(\frac{2\sqrt{n}}{\sqrt[4]{p}}-\xi\right)+\| \beta^* \|_1\left(1-\frac{\sqrt{n}}{\sqrt[4]{p}}\right).
\]

Considering $\Delta_{\xi}=\left\{u\in\mathbb{R}^p, \|{u^*}^c\|_1\xi\leq\|{u^*}\|_1\left(\frac{2\sqrt{n}}{\sqrt[4]{p}}-\xi\right)+\| \beta^* \|_1\left(1-\frac{\sqrt{n}}{\sqrt[4]{p}}\right)\right\}$ and assuming that \\ $k_{\xi}^* = \min\limits_{u\in \Delta_{\xi}}\frac{\frac{\sqrt{p^*}}{\sqrt{n}}\|u\|_X}{\frac{2\sqrt{n}}{\xi\sqrt[4]{p}}\|{u^*}\|_1+\frac{\| \beta^* \|_1}{\xi}\left(1-\frac{\sqrt{n}}{\sqrt[4]{p}}\right)}>k>0$ we get

\begin{eqnarray*} 
L(\hat{\beta}(\lambda))-L(\beta^*) & \le & | L(\hat{\beta}(\lambda))-L(\beta^*) | \le \frac{\lambda\sqrt{n}}{c\sqrt[4]{p}}(\|{u^*}\|_1+\|{u^*}^c\|_1)\\
& \le & \frac{\lambda \sqrt{n}\|u\|_1 } {c\sqrt[4]{p}} \le \frac{\lambda \sqrt{p^*} \| u \|_X }{ \sqrt[4]{p} k_{\xi}^* } .
\end{eqnarray*}
 The rest of the proof is virtually identical with the last part of the argument we detailed for 
 Theorem \ref{oracle}.

\end{proof}
}

\end{document}